\documentclass{article} 
\usepackage{enumerate}  
\usepackage{a4}
\usepackage{times}
\usepackage{amsmath}
\usepackage{array}
\usepackage{amssymb}
\usepackage{graphicx}
\usepackage[left=2.75cm,right=2.25cm,top=2.9cm,bottom=2.25cm]{geometry}

\expandafter\chardef\csname pre amssym.def  at\endcsname=\the\catcode`\@
\catcode`\@=11
 
\def\undefine#1{\let#1\undefined} 
\def\newsymbol#1#2#3#4#5{\let\next@\relax
 \ifnum#2=\@ne\let\next@\msafam@\else
 \ifnum#2=\tw@\let\next@\msbfam@\fi\fi
 \mathchardef#1="#3\next@#4#5}
\def\mathhexbox@#1#2#3{\relax
 \ifmmode\mathpalette{}{\m@th\mathchar"#1#2#3}%
 \else\leavevmode\hbox{$\m@th\mathchar"#1#2#3$}\fi}
\def\hexnumber@#1{\ifcase#1 0\or 1\or 2\or 3\or 4\or 5\or 6\or 7\or 8\or
 9\or A\or B\or C\or D\or E\or F\fi}
 
\font\tenmsa=msam10
\font\sevenmsa=msam7  
\font\fivemsa=msam5
\newfam\msafam
\textfont\msafam=\tenmsa
\scriptfont\msafam=\sevenmsa
\scriptscriptfont\msafam=\fivemsa 
\edef\msafam@{\hexnumber@\msafam}

\font\tenmsb=msbm10
\font\sevenmsb=msbm7
\font\fivemsb=msbm5
\newfam\msbfam
\textfont\msbfam=\tenmsb
\scriptfont\msbfam=\sevenmsb
\scriptscriptfont\msbfam=\fivemsb
\edef\msbfam@{\hexnumber@\msbfam}

   \font\tengothic=eufm10
   \font\sevengothic=eufm7
   \newfam\gothicfam
   \textfont\gothicfam=\tengothic
   \scriptfont\gothicfam=\sevengothic
   \def\goth#1{{\fam\gothicfam #1}}
   \font\tenmsb=msbm10
   \font\sevenmsb=msbm7
   \newfam\msbfam
   \textfont\msbfam=\tenmsb
   \scriptfont\msbfam=\sevenmsb
   						
\newtheorem{prop}{Proposition}[section] 
\newtheorem{set}[prop]{Setting} 
\newtheorem{rem}[prop]{Remark}
\newtheorem{thm}[prop]{Theorem}
\newtheorem{coro}[prop]{Corollary}

\newtheorem{lemma}[prop]{Lemma}
\newtheorem{notat}[prop]{Notation}

\newtheorem{(*)}[prop]{}

\newcommand{\af}{\mbox{\sf  {\it I}\hspace{-0.65 mm}A}}

\newcommand{\nat}{\mbox{${\rm I\hspace{-.6 mm}N}$}}

\newcommand{\p}{\mbox{$\  \longrightarrow \ $}}

\newcommand{\enu} {\begin{enumerate}[{\rm(1)}]} 
\newcommand{\enua} {\begin{enumerate}[ $(a)$]} 
\newcommand{\enui} {\begin{enumerate}[ $(i)$]} 
\newcommand{\denu} {\end{enumerate}}

\font\tengothic=eufm10 
\font\sevengothic=eufm7
\newfam\gothicfam 
\textfont\gothicfam=\tengothic
\scriptfont\gothicfam=\sevengothic
\def\goth#1{{\fam\gothicfam #1}}

\begin{document}
\begin{center}

\renewcommand{\thefootnote}{\fnsymbol{footnote}}
\addtocounter{footnote}{0}

{\LARGE  Explicit minimal resolution for certain monomial curves.}
\end{center}

\begin{center}

\renewcommand{\thefootnote}{\fnsymbol{footnote}}
\addtocounter{footnote}{0}

{\large
    Anna Oneto,   Grazia~Tamone . 
}  \\[2mm]
{\small DIMA - University of Genova,
via Dodecaneso 35, I-16146 Genova,  Italy.}\\
{\small E-mail:    oneto@diptem.unige.it,   tamone@dima.unige.it }
\end{center} 

\vspace*{-5mm}

\noindent \rule[0pt]{\textwidth}{1pt}

\vspace*{-3mm}

{\baselineskip8pt

\begin{abstract}

{\baselineskip5pt

\noindent With a view to study problems of smoothability , we construct a minimal free resolution for the coordinate ring of an algebroid monomial curve associated to an $AS$ numerical semigroup (i.e. generated by an arithmetic sequence), obtained independently of the result of  \cite{gss2} and  equipped with the explicit description of all the involved maps. }
\end{abstract}
}

\vspace*{-2mm}

\noindent \rule[0pt]{\textwidth}{1pt}

{\baselineskip6pt
\vskip4pt

\noindent {\footnotesize\it Keywords~}{\scriptsize\sf :} {\scriptsize\sf
 Numerical semigroup, Arithmetic sequence, Monomial curve, Free resolution, Betti number. }
\smallskip

\noindent {\footnotesize\it {\rm 1991} Mathematics Subject Classification~}: {\scriptsize\sf 
13D02.}

}

\vspace*{-5mm}

\renewcommand{\thefootnote}{\arabic{footnote}}
\addtocounter{footnote}{0}

\setcounter{section}{-1}

 \section{Introduction}
  Let $k [x_0,...,x_n]/I\simeq k[t^s,\,\, s\in S]  $ be the affine coordinate ring of a monomial algebraic curve $X\subseteq\af^{n+1}_k$ defined by a numerical semigroup $S$ and let $R=k[[x_0,...,x_n]]/I$ \, ($k$ field).  Several properties of the associated graded ring ${\cal G}$ of $R$ have been studied in the recent years; on the other side, some questions related to the homological invariants of the local ring $R$ are still open. When $S$ is generated by an arithmetic sequence $m_0,...,m_n$ ($AS$ semigroup) the generating ideal $I$ of $X$ has a  nice shape: in fact it is generated by the $2\times 2$ minors of two matrices, as first noted in \cite{gss1} and  the Betti numbers of ${\cal G}$ are calculated in   \cite{lr}.\ \
\par  In the recent paper \cite{gss2} the authors give a minimal free resolution of the ring $k [x_0,...,x_n]/I$,   based on this \lq\lq bideterminantal'' shape of the ideal $I$ and the mapping cone procedure. 
In particular they deduce the Betti numbers of $I$   \cite[Theorem 4.1]{gss2}.
 Essentially with an analogous technique, but independently,  we have reached the same result: the difference is our more explicit definition of the maps involved in this construction.
 This  punctual description is a basic tool to obtain a determinantal characterization of the first syzygies of the ideal $I$, which turns out to be useful in the study of the smoothability of these monomial curves (work in progress).
 \section{Setting.}
 \begin{notat} \label{nota0} 
{\rm \enu\item \  For a semigroup $S$ generated by an arithmetic sequence,  $S=\sum_{0\leq i\leq n }\nat m_i$ , where $m_i=m_0+i\,\! d$,   $(1\leq i\leq n$) \, and  \, $GCD(m_0,d)=1$,  let $\,  a ,b,\mu  \in {\mathbb  N}$ \ \ be such that \ \ \\
\centerline{ $m_0=an+b$ , \ \quad$a\geq 1$,\quad$1\leq b\leq n,\quad\mu:=a+d$.\hspace{1cm}}\\
 Let \, $P:=k[x_0,...,x_n]$ \, ($k$ field), \, with \, $weight(x_i):=m_i$, \, and let $k[S]=k[t^s,\,s\in S]$.\\  The     defining ideal $I\subseteq P$ of the curve  $X=Spec (k[S])$ (shortly  {\it AS monomial curve} ) is generated by the $2\times 2$ minors of the following two matrices: \\ [2mm] 
   \centerline{ $A:=\left(\begin{array}{ccccc}
 x_0&x_1 \ \dots \     x_{n-2}&x_{n-1}\\
 x_1&x_2 \  \dots\     x_{n-1}&x_{n }\\
 \end{array}\right)\quad
A':=\left(\begin{array}{ccccc}
 x_n^a&x_0\ \dots\  x_{n-b}\\
 x_0^{\mu}&x_b \ \dots\  \ x_{n}  \ \ \\
 \end{array}\right)    $.}\\[2mm]
  and   a minimal set of generators for $I$ can be obtained by  the $\displaystyle{\binom{n}{2}}$ 
 maximal minors $\, \{ f_1,...,f_{\binom{n}{2}}\}$ \, (we choose   lexicographic order) of the matrix $A$ and the $n-b+1$ maximal minors $M_{1   j}$  containing the first column of the matrix $A'$ \, $($see {\rm \cite[Theorem 1.1]{gss1}\ )}. 
 \item
  We call  $  {\goth C} $ the codimension two ideal generated by the $2\times 2$ minors of the matrix $A$ (which is   the ideal of the cone over the rational normal curve of ${\mathbb P}^n$).\ 
 \item
For $h=0,...,n-b$,   we shall denote by $\,g_h\,$ the minor  $det \left(\begin{array}{llccc}
 x_n^a&   x_{n-b-h}\\
 x_0^{\mu}&  x_{n-h}  \ \ \\
 \end{array}\right) $ of $A'$    and by $\delta_h$ its $weight$: \\
 {  $g_0:= x_n^{a+1}-x_0^{\mu} x_{n-b},\quad g_h= x_n^ax_{n-h}-x_0^{\mu} x_{n-b-h}, \quad g_{n-b}:=   x_n^{a}x_b -x_0^{\mu+1},\quad \delta_h=am_n+m_{n-h} $.}

 \denu }\end{notat}
 \begin{lemma}\label{x0}
 With { Notation \ref{nota0}}: \quad $ \big(  {\goth C} +(g_0,...,g_h)\big):g_{h+1}=(x_1,...,x_n) ,\quad  $ for each $  \, h=0,...,n -b-1 $.
  \end{lemma}
 Proof. First observe that  
  \ \ $x_{i+1}g_{h+1}-x_ig_h= x_{i+1}(x_n^ax_{n-h-1}-x_0^{\mu}  x_{n-b-h-1} )-x_i(x_n^ax_{n-h}-x_0^{\mu} x_{n-b-h}) =x_n^a(x_{i+1} x_{n-h-1}-
 x_{i} x_{n-h})-x_0^{\mu}(x_{i+1} x_{n-b-h-1}-x_i x_{n-b-h})\in   {\goth C} $ for each $i=0,...,n-1,\ \ h=0,...,n-b-1$.
Hence the inclusion $\supseteq$ is clear.\\
Now assume that $x_0 g_{h+1}\in   {\goth C} +(g_0,...,g_h)$. Then $x_0 g_{h+1}=\beta+\alpha_0 g_0+...\alpha_h g_h,  \alpha_i\in P, \beta\in   {\goth C} . $
Hence \\[1mm] $x_0^2  g_{h+1}=x_0\beta+\alpha_0x_0 g_0+...\alpha_h x_0 g_h= \beta_1+\alpha_0x_1 g_1+...+\alpha_h x_1 g_{h+1}=\beta_2+\alpha_0x_2 g_2+...+(\alpha_{h-1}x_2+\alpha_h x_1 )g_{h+1}=...=\beta_h+ \alpha g_{h+1}$, with
$\beta_h\in   {\goth C} ,\ \ \alpha\in (x_1,...,x_n).$ This would imply that $(x_0^2-\alpha)g_{h+1}\in   {\goth C} $, impossible since $  {\goth C} $ is prime, $g_{h+1}\notin   {\goth C} ,\, x_0^2-\alpha=x_0^2-(\alpha_0x_{h+1}+\alpha_1x_{h }+...+\alpha_h x_{ 1})\notin   {\goth C}  $  (because $t^{2m_0}\notin (t^{m_1},...,t^{m_n}$)\,). \quad$\diamond$

\section { Free resolution of the ideal  $I$.}

By means of the mapping cone technique, starting with the Eagon-Northcott  resolution ${\mathbb E} $ of the ideal ${\goth C}$ and the Koszul complex ${\mathbb K}$ of $P/(x_1,...,x_n)$, we can construct a free (non-minimal in general) resolution of the ideal $I$. This resolution is a lifting of the one found in \cite{lr} for the associated graded ring ${\cal G}$ of the curve $X$: in particular we'll see that the Betti numbers of ${\cal G}$ and    $A$ are the same. We recall the main tools:
 \begin{set}\label{nota1}{\rm With Notation \ref{nota0},  let $R=k[[x_0,...,x_n]]/I$. 
 \enu
  \item   The Eagon-Northcott   free resolution    for the $R$-ideal $  {\goth C} =\left(f_1,...,f_{\binom{n}{2}}\right)$,  is the complex 
  \\ \centerline{ $\begin{array}{lccccccccccccccccc}
&& &&d_{s} {}& &&&d_2 {}&&d_1 {}&& \\[-0.5mm]   
 {\mathbb E}:&0\p E_{n-1}\p &...\!\!\longrightarrow\!\!&  E_{s}    &\!\!\longrightarrow\!\!& E_{s-1}  &... & E_{2}  &\!\!\longrightarrow\!\!& E_{1}  &\!\!\longrightarrow\!\!& E_0&\!\!\longrightarrow\!\!& P/  {\goth C} &\!\!\longrightarrow\!\!& 0 \
    \end{array} $}\\[2mm]
  where $E_0\simeq P, \quad\,E_{s}=\wedge^{s+1} P^n\bigotimes \left( Sym_{s-1}(P^2)\right)^*\simeq P^{\beta_s}(-s-1),$  \,  for \,$1\leq s\leq n-1,$ \quad   $ \beta_s=s\binom{n}{s+1}$,  \\[2mm] $Sym_{s-1}(P^2)$ is a free $P$-module of  $rank \ s$ and basis $\{\lambda_0^{v_0 }\lambda_1^{v_1 }| \, v_0+v_1=s-1\}$ (  $\, 1\leq s\leq n-1)$   (see \cite{lr}).\\[2mm]
Let \quad $<e_{i_1}{}\wedge...\wedge e_{i_{s+1}}{}\otimes \lambda_0^{v_0}\lambda_1^{v_1}, \ ( 1\leq i_1<i_2<...<i_{s+1}\leq n,\quad v_0+v_1=s-1 )>$  be   the basis of $E_s$ and let $\varepsilon$ be the basis of $E_0$.\quad   
The maps in $\mathbb E$ are:\\[2mm]
\indent $d_1{}:E_1\p E_0,$  \, $  e_{i_1}{}\wedge e_{i_2}{}\mapsto (x_{i_1-1}x_{i_2}-x_{i_1}x_{i_2-1})\varepsilon,\quad 1\leq i_1<i_2\leq n$,  \\[2mm]
      \indent $d_s{}\left( (  e_{i_1}{}\wedge...\wedge  e_{i_{s+1}}{}) \otimes \lambda_0^{v_0  }\lambda_1^{v_1 }
\right)=\Delta_0( e_{i_1}{}\wedge...\wedge e_{i_{s+1}}{})\otimes \lambda_0^{v_0-1 }\lambda_1^{v_1 }+\Delta_1( e_{i_1}{}\wedge...\wedge e_{i_{s+1}}{}) \otimes \lambda_0^{v_0 }\lambda_1^{v_1-1 },\quad s\geq 2 ,$\\[2mm]
\noindent   where      only the summands with non-negative powers of $\lambda_0,\lambda_1$ are considered, and, \,   for \, $q=0,1, s\geq 1$, the  maps\\ 
$\Delta_q:\wedge^{s } P^n\p\wedge^{s-1} P^n$ are defined as: \quad\\[2mm]
 \centerline{ $\begin{array}{l}\Delta_q(  e_{i_1}{}\wedge...\wedge e_{i_{s}}{} ) :=\sum_{j=1}^{s}(-1)^{j+1} x_{i_{j}+q-1}   \, e_{i_1}{}\wedge...{\widehat {e_{i_j}{}}}...\wedge e_{i_{s}}{} ,\quad (1\leq i_1<...<i_s\leq n ,\quad  2\leq s\leq n),\\
 \Delta_q(  e_{i } ) := x_{i+q-1}\varepsilon,  \quad (1\leq i \leq n ,\quad s=1).\end{array}$.} 
 \item  The   Koszul  complex  ${\mathbb K}$,   minimal free resolution for $P/(x_1,...,x_n)$, is:\\ 
 {\centerline{$\begin{array}{lccccccccccccccccccccc}
 &&_{\Delta_1}&&&&_{\Delta_1}\\[-.3mm]
{\mathbb K}:&  0 \p   K_{n}  &\p&...  &\p &K_{1}&\p   &K_0&\!\!\longrightarrow\!\!& P/(x_1,...,x_n) \p  0 
 \end{array} $}}
 \item For $ \,1\leq h\leq n-b)$  consider tha ideal $  {\goth C} _{h-1}:=(  {\goth C} ,g_0,g_1,...,g_{h-1})\subseteq P$. \quad     By (\ref{x0}), \, $P/(  {\goth C} _{h-1}:g_h)=P/(x_1,...,x_n)$;
 hence the Koszul complex      gives a minimal free resolution $\,\,  \mathbb K(-\delta_h)$ \, of  \, $P/(x_1,...,x_n)$:\\ 
{\centerline{ $\begin{array}{lccccccccccccccccccccc}
  &&& &d'_{s}& && &\cdot\, g_h\\[-2mm]  
 {\mathbb K}(-\delta_h) :&  0\p   K_{n}(-\delta_h)   &...   &K_{s}(-\delta_h)&\p &...  &K_0(-\delta_h)&&\longrightarrow\!\!& P(-\delta_1)/(  {\goth C} _{h-1}:g_h) \p  0 
 \end{array} $}}\\[2mm]
     where $\,\ d'_s  =\Delta_1  $,  $\quad (1\leq s\leq n)$. 
\item  For $( \,0\leq h\leq n-b)$ denote respectively by \\[2mm] 
 \centerline{  $ \begin{array}{cll} 
\varepsilon_h, &{\rm  the \,  (canonical) \, basis\, of}& K_0(-\delta_h),   \\[2mm] 
 e^{(h)}_{i_{1}} \wedge...\wedge e^{(h)}_{i_{s}},   &{\rm  the \,  (canonical) \, basis\, of}& K_{s}(-\delta_h)  \quad (1\leq h\leq n-b),   \\[2mm]   
 \varepsilon_0,  &{\rm  the \,  (canonical) \, basis\, of}&  E_0(-\delta_0),\\[2mm]
e^{(0)}_{i_{1}} \wedge...\wedge e^{(0)}_{i_{s}}  &{\rm  the \,  (canonical) \, basis\, of}&E_{s-1}(-\delta_0)\quad (1\leq i_1<...<i_s\leq n  )   
\end{array}$.} 
 \item[] Further, for simplicity, when no confusion on indices occurs, we shall write \\[2mm] 
   \centerline{$ W:= e_{i_{1}}\wedge...\wedge e_{i_{s}},    \quad  \quad  W^{(h)}=e^{(h)}_{i_1}\wedge...\wedge e^{(h)}_{i_{s}} , \quad (0\leq h\leq n-b)  $.\quad} 
\item   Recalling that \, $weight(x_i)=m_i=m_0+i\,d \,\, (0\leq i\leq n)$, \, set \quad $ weight (\lambda_0):=0,\quad weight (\lambda_1):=d,$ \\[2mm]
\centerline{$  weight (\varepsilon_h)=-\delta_h, \,\, weight (e_{i_{1}}^{(\,)}\wedge...\wedge e_{i_{s}}^{(\,)})=m_{i_{1}}+...+m_{i_{s}}-(s-1)d. $  }\\[2mm]
 Therefore the modules $K_s, E_s$ are graded as follows:
 
\item[ ]
 $K_s(-\delta_h)=\bigoplus_{1\leq i_1<...<i_s\leq n}P\left(-\delta_h-m_{i_1 }-...-m_{i_s }+(s-1)d\right)$, \, for \quad $1\leq s\leq n+1,\quad h=0,...,n-b,$
 
 \item[ ] $E_s =\bigoplus_{0\leq v_1\leq s-1}\left[\bigoplus_{1\leq i_1<...<i_{s+1}\leq n}P\left( -m_{i_1 }-...-m_{i_{s+1} }+(s-v_1)d\right)\right]$, \,\, for \quad$1\leq s\leq n-1,$\quad
   $E_0\simeq P$
    {\rm
\item With the preceding assumptions we  shall  define by $:$ \\[2mm]
  \centerline{$\begin{array}{llll}
  F_0^{(0)}&:=&E_0 \\
F_s^{(0)}&:=& E_{s-1}(-\delta_0)\oplus E_{s }  \quad  (1\leq s\leq n-1)\\ 
F_n^{(0)}&:=& E_{n-1}(-\delta_0)\\
  F_0^{(h)}&:=&E_0(-\delta_h)\quad (h\geq 1)\\
    F_s^{(h )}&:=&K_{s-1}(-\delta_{h })\oplus...\oplus K_{s-1}(-\delta_{1})\oplus E_{s-1}(-\delta_0)\oplus E_s \quad (s,h\geq 1)\\
  \end{array}$.}
 \\[2mm] 
 $\psi_{s}^{(0)}:E_s(-\delta_0)\p E_s $, \quad the multiplication by $g_0$,\,\, for all $s\geq 0$;\\[2mm] 
  $ \psi_0^{(h)}:K_0(-\delta_h)\p F_0^{(h-1)}$   the multiplication by $\,\, g_h=x_n^ax_{n-h}-x_0^{\mu}x_{n-b-h};\quad (h\geq 1) $  
 \\[2mm]  
 $  \psi_{s}^{(h)}:K_s(-\delta_h)\p K_{s-1}(-\delta_{h-1})\!\oplus\!  E_{s-1}(-\delta_0)\!\oplus\! E_s\subseteq F_s^{(h-1)}\quad(s,h\geq 1) $:
 \enu
 
 \item[] $\psi_1^{(1)}(e_i^{(1)}) =x_{i-1}\varepsilon_{0}\oplus e_i\wedge \left( x_0^{\mu}e_{n-b}  -x_n^ae_{n} \right) \in   E_{0}(-\delta_0)\!\oplus\! E_1$,\quad $i=1,...,n , \qquad (h=1)$;
  \item[]$\psi^{(h)}_{1}(e^{(h)}_{i } )=\ \
    -\Delta_0(e^{(h-1)}_{i } ) + \phi^{(h)}_1 (e^{(1)}_{i } ),\qquad (h\geq 2)    $;

  \item[]  $\psi_{s}^{(1)}(W ^{(1)} )= $
$-  W^{(0)}\otimes \lambda_1^{s-2} \oplus    \phi_s^{(1)}(W^{(1)}),\qquad (s\geq 2);$ 
 \item[] $ \psi^{(h)}_{s}(W^{(h)})=     (-1)^{s-1}\Delta_0(W^{(h-1)})+(-1)^h W^{(0)}\otimes\lambda_0^{h-1}\lambda_1^{s-h-1}+ \phi^{(h)}_s (W^{(h)}  ),\qquad (s,h\geq 2)  $\denu
 \indent {\rm where} \,$\phi^{(h)}_s :K_s(-\delta_h)\p E_s$ \, is  \\[2mm] \indent
 $\phi^{(h)}_s(W^{(h)} ) :=W\wedge\sum_{k=1}^h(-1)^k   ( x_n^a e_{n+k-h} ^{}-x_0^{\mu} e_{n-b+k-h} ^{})\otimes \lambda_0^{k-1}\lambda_1^{s-k} ; $
 
  \item     The maps \quad$d_s^{(h)}: F_{s}^{(h)}\p F_{s-1}^{(h)}$ are defined as follows: 
    \item[]
 $d_1^{(0)}(\alpha\varepsilon_{0}+\beta  e_i\wedge e_j) =[\alpha g_0+\beta (x_{i-1}x_{j}-x_{i}x_{j-1})]\varepsilon$;\
 \item[] For $s=1, h\geq 1$, with
    $\begin{array}{cccccccc} F_1^{(h )}=K_{0}(-\delta_{h })\oplus...\oplus K_{0}(-\delta_{1})\oplus E_{0}(-\delta_0)\oplus E_1, \end{array}  \quad   
 F_0^{(h )}=  E_0(-\delta_{h}) , $ \\
 $ d_{1 }^{(h)}: F_{1}^{(h)}\p F_{0}^{(h )} $  \,
 is  the component-wise product 
 \item[] \centerline{ $(\alpha_h \varepsilon_h\oplus ...\oplus \alpha_0 \varepsilon_0\oplus \beta\, e_i\wedge e_j)\cdot\left(
 g_h\oplus    g_{h-1}\oplus...   \oplus   g_0\oplus (x_{i\!-\!1}x_{j}\!-\! x_{j\!-\!1}x_{i}) 
\right)$;}  
 $d_{2 }^{(0)}\left((e_i\wedge e_j)\oplus (e_h\wedge e_k\wedge e_l \otimes  \lambda_q )\right)\!=\! 
(x_{i\!-\!1}x_{j}\!-\!x_{i}x_{j\!-\!1})\varepsilon_{0}\oplus\left( -g_0 \ e_i\wedge e_j +d_2(e_h\wedge e_k\wedge e_l\otimes \lambda_q )\right),$ $\,(q=1,2) $;  
\item[] $d_{s }^{(0)}:E_{s-1}(-\delta_0)\oplus E_{s}\p E_{s-2}(-\delta_0)\oplus E_{s-1},\quad s\geq 2,\quad$ {\rm is  } \ \
$      \left(
\begin{array}{ccccccc}
  d^{ }_{s-1}&0       \\
  (-1)^{s-1}(_{...} \cdot\, g_0)& d^{ }_{s}  \\
     \end{array}
\right);$

\item[]
$d_{s }^{(1)}:K_{s-1}(-\delta_1)\oplus E_{s-1}(-\delta_0)\oplus E_{s}\p K_{s-2}(-\delta_1)\oplus E_{s-2}(-\delta_0)\oplus E_{s-1},$\quad $(s\geq 2)$  is given by
\item[] $   \left(
\begin{array}{ccc}
  d'_{s-1}&0       \\
  (-1)^{s-1}\psi^{(1)}_{s-1}& d^{(0)}_{s}  \\
     \end{array}
\right)\equiv  \left(
\begin{array}{ccc}
  \Delta_{1,s-1}&0   &0    \\
  -(_{...}\otimes \lambda_1^{s-3})&  d_{s-1} &0       \\
     \\
   (_{...}\wedge \left(x_0^{\mu}e_{n-b} 
-x_n^a  e_ {n} \right)\otimes \lambda_1^{s-2})& (-1)^{s-1}(_{...}\cdot\, g_0)& d_{s}^{}  
        \end{array}\right)$

\item[] $ d_{s }^{(h)}: F_{s}^{(h)}\p F_{s-1}^{(h )} $,\quad is associated to the matrix   \ \
 $ 
 \left(
\begin{array}{ccc}
  d'_{s-1}&0       \\
  (-1)^{s-1}\psi^{(h)}_{s-1}& d^{(h-1)}_{s}  \\
     \end{array}
\right),$ $\quad (2\leq s\leq n, \quad h\geq 1)$  
\\
 $ 
 d_{s }^{(h)}:  \left(
\begin{array}{cccccc}
 \Delta_{1 }&0   &...&0&0  &0   \\

(-1)^{s-1}\psi_{s-1}^{(h)}&  \Delta_{1 }&...&0&0   &0    \\
0&(-1)^{s-1}\psi^{(h-1)}_{s-1}&...&0&0&0\\
...&&...  &0&0& 0 \\
...& &...  &\Delta_1&0&0 \\
 0& 0&...&(-1)^{s-1}\psi_{s-1}^{(1)}&  d_{s-1}^{} &0      \\

     0& 0&...&0& (-1)^{s-1}(\,\,\, \cdot g_0)& d_{s}^{}  
     \end{array}\right)$. 
   \item[] Hence, \, for \quad $1\leq k\leq h$,\ \\[2mm]
     $\begin{array}{lll}d_{s+1}^{(h)}(W^{(k)})  &= \Delta_1(W^{(k )})\oplus (-1)^{s}\psi_{s}^{(k)}(W^{(k)}) \in K_{s-1}(-\delta_k) \oplus K_{s-1}(-\delta_{k-1})\oplus  E_{s-1}(-\delta_0)\oplus E_s ;\\
 &= \Delta_1(W^{(k)})\oplus (-1)\Delta_0(W^{(k-1)})\oplus (-1)^{s+k}W^{(0)}\otimes \lambda_0^{h-1}\lambda_1^{s-h-1}\oplus
      \phi^{(h)}_s (W  )\end{array}\\
    \begin{array}{lll}     d_s^{(h)}(W^{(0)}\otimes \lambda_0^{v_0} \lambda_1^{v_1})&=& \Delta_0(W)\otimes \lambda_0^{v_0-1} \lambda_1^{v_1}+\Delta_1(W)\otimes \lambda_0^{v_0} \lambda_1^{v_1-1}+(-1)^{s-1 }W^{(0)}\otimes \lambda_0^{v_0} \lambda_1^{v_1}g_0,\\   &&(v_0+v_1=s-2);\\  d_{s+1}^{(h)}(W \otimes \lambda_0^{v_0} \lambda_1^{v_1})&=& d_{s+1} (W \otimes \lambda_0^{v_0} \lambda_1^{v_1}),\quad (v_0+v_1=s-1).
 \end{array}$} 
 \denu 

}  \end{set}

\subsection{Mapping cone construction} 
Now we   construct a   free resolution of the ideal $I$.
\begin{prop}\label{cs}
With notation {\rm (\ref{nota0}), (\ref{nota1}}), \, the  following complex is a free resolution of the ideal $I$:\\[2mm] 
\centerline{$\begin{array}{cccccccccccccccccccccc}
& &&  d_{s}^{(n-b)}  & &&&d_{2}^{(n-b)} &&d_{1}^{(n-b)}  &&  \\[-1.5mm] 
 {\mathbb F}_{n-b} :&&  0 \p    C _{n+1}  \p ...\p C_{s} &\p& C_{s-1}   &...  & C _2  &\!\!\longrightarrow\!\!& C _1  &\!\!\longrightarrow\!\!& P   \end{array} $}\\[2mm]
where $C_s=F^{(n-b)}_{s}=  K_{s-1}(-\delta_{n-b })\oplus...\oplus K_{s-1}(-\delta_{1})\oplus E_{s-1}(-\delta_0)\oplus E_s  $, \\[2mm]
$ dim(C_s)= (n-b)\binom{n}{s-1}+(s-1)\binom{n}{s}+ s\binom{n}{s+1}.$
\end{prop}
The proof consists of several steps. \\[2mm] 
{\bf STEP 0.} (Well-known). Let $g_0:= x_n^{a+1}-x_0^{\mu} x_{n-b},\quad \delta_0=\, deg(g_0)$. The following diagram is commutative:\\[1mm]
 \centerline{$\begin{array}{cccccccccccccccccc}
 & &&&&&&&&0\\
 &  &&&&&&&&\downarrow \\
&&&&d^{}_{2} & &d^{}_1&&d^{}_0\\[-1mm]  
 {\mathbb E} (-\delta_0) :  &   E_{s}(-\delta_0)& ... & E_{2}(-\delta_0)  & \!\!\longrightarrow\!\!& E_{1}(-\delta_0)&\!\!\longrightarrow\!\!&E_0(-\delta_0)&\!\!\longrightarrow\!\!& P(-\delta_0)/  {\goth C} &\!\!\longrightarrow\!\!& 0\\
\\
  &\ \ \ \ \downarrow\psi_{s }^{(0)}&&\ \ \ \ \downarrow\psi_{2}^{(0)}&&\ \ \ \ \downarrow\psi_1^{(0)}&&\ \ \ \ \downarrow\psi_{0}^{(0)}&&\ \ \ \ \downarrow \cdot g_0 \\
 &&&&d_{2} ^{}&&d_1 ^{}&&d_0 ^{}\\[-1mm]  
 {\mathbb E} : &  E_{s}    & ...& E_{2}    &\!\!\longrightarrow\!\!& E_{1}  &\!\!\longrightarrow\!\!& E_0&\!\!\longrightarrow\!\!& P/  {\goth C} &\!\!\longrightarrow\!\!& 0\\
 \\
 &&&&&& &&\searrow&\downarrow \\ \\  
  &&&&&&&&&P/(  {\goth C} ,g_0) \\  
  && &&&&&&&\downarrow \\
&&&&&&&&&0
  \end{array} $}
     where the $\psi_i^{(0)}$ are the multiplication by $g_0$.  By mapping cone one constructs the complex:
 $$\begin{array}{llcccccccccccccccc}
  &\qquad\qquad \qquad\qquad\qquad \qquad  \qquad\qquad\qquad\qquad d_{s}^{(0)}&   \ \ &&   &\ \ \ \\[-1mm]

 {\mathbb F}_0:& 0 \p   E_{n-1}(-\delta_0) \p...\p E_{s-1}(-\delta_0)\oplus E_{s } \p   E_{s-2}(-\delta_0)\oplus E_{s-1}\  & ...& E_0(-\delta_0)\oplus E_1\p E_0\\
&\qquad\qquad\|\qquad\qquad\qquad\qquad\qquad\qquad\quad\| &&\!\!\|&& \  \\
&\qquad\qquad F_n^{(0)} \qquad\qquad\qquad\qquad\qquad\quad   F_s^{(0)}&  &F_1^{(0)} \\  \end{array} $$
   $d_1^{(0)}(\alpha\varepsilon_0 +\beta e_i\wedge e_j) =[\alpha g_0+\beta(x_{i-1}x_{j}-x_{i}x_{j-1})]\varepsilon$,\quad  \\[2mm]  
  $d_{s }^{(0)}:E_{s-1}(-\delta_0)\oplus E_{s}\p E_{s-2}(-\delta_0)\oplus E_{s-1},\quad s\geq 2,\quad$  defined by 
$   \left(
\begin{array}{ccc}
  d^{ }_{s-1}&0       \\
  (-1)^{s-1}\psi^{(0)}_{s-1}& d^{ }_{s}  \\
     \end{array}
\right)$.\quad   \\[2mm]
{\bf STEP 1.}  Let $  {\goth C} _0=(  {\goth C} ,g_0)\subseteq P$ and let  $g_1= x_n^ax_{n-1}-x_0^{\mu} x_{n-b- 1}$. Consider the Koszul complex   $\mathbb K(-\delta_1)$ of   $P(-\delta_1)/(  {\goth C} _0:g_1)$; 
  from the exact sequence\quad$0\p P/((  {\goth C} ,g_0):g_1) \p P/(  {\goth C} ,g_0)\p P/(  {\goth C} ,g_0,g_1)\p 0$, \,\,  we obtain  the commutative diagram \\ 
\centerline{ $\begin{array}{rrcccccccccccccccccccc}
 && &&&& &&0\\
 &&  & &&&&&\downarrow \\
 &         &  &&d'_s&&&d'_0\\  
 {\mathbb K}(-\delta_1) :&   0\p   K_{n}(-\delta_1)&...  &  K_{s}(-\delta_1)&\p&... \p&K_0(-\delta_1)&\!\!\longrightarrow\!\!& P(-\delta_1)/(  {\goth C} _0:g_1)&\!\!\longrightarrow\!\!& 0\\
\\
 &\downarrow\psi_{n }^{(1)}&& \downarrow\psi_s^{(1)}&&&\downarrow\psi_{0}^{(1)}&&\ \ \ \ \downarrow \cdot g_1 \\
 &&  && d_s ^{(0)}&& &d_0 ^{(0)}  \\ 
{\mathbb F}_0:&  0\p   F^{(0)}_{n}  &...& F^{(0)}_{s}&\p&... \p& F^{(0)}_0\ \ \   \ \ \ &\!\!\longrightarrow\!\!& P/  {\goth C} _0&\!\!\longrightarrow\!\!& 0\\
\\ &\|\ \ \ \ \ \ \ \ \ &&\|&&&   \|\ \ \ \ \ \ \ \ \ &\searrow&\downarrow \\ \\
  &E_{n-1}(-\delta_0) &&E_{s-1}(-\delta_0)\oplus E_s&&&E_0\ \ \ \ \ \ \ \ &&P/(  {\goth C} ,g_0,g_1) \\ \\
   &&&&&&&&\downarrow \\
 &&&&&&&&0
\end{array} $}
  The maps, as defined  in (\ref{nota1}.6),  are :
\enu

 \item[] $ \psi_0^{(1)}$, \quad   the multiplication by $\,\, g_1=x_n^ax_{n-1}-x_0^{\mu}x_{n-b-1},$ 
 \item[] $\psi_1^{(1)}(e_i^{(1)}) =x_{i-1}\varepsilon_0 +e_i\wedge \left( x_0^{\mu}e_{n-b}  -x_n^ae_{n} \right) $,\quad $i=1,...,n$, 

\item[]$\psi_{s}^{(1)}(e_i^{(1)} )= $
$ -e_i^{(0)}\otimes \lambda_1^{s-2}\oplus e_i \wedge \left( x_0^{\mu}e_{n-b} -  x_n^a   e_ {n}  \right)\otimes \lambda_1^{s-1},\quad s\geq 2.$

\denu
One can check directly the commutativity of the diagrams for each $(s=1,...,n)$
$$  \begin{array}{cccccccccccccccccccccc}

 & d'_{s}&&\\  
    K_{s}(-\delta_1)&\!\!\longrightarrow\!\!& K_{s-1}(-\delta_1)\\
\\
\ \ \ \ \  \downarrow\psi_{s }^{(1)}& &\downarrow\psi_{s-1 }^{(1)} \\

  &  d_{s} ^{(0)}   \\ 

    E_{s-1}(-\delta_0)\oplus E_s             &\!\!\longrightarrow\!\!& E_{s-2}(-\delta_0)\oplus E_{s-1}\\\

  \end{array} $$ 
 If $s=1$,\quad $\psi_{0 }^{(1)} d'_{1}(e_i^{(1)})=x_i g_1\varepsilon$,\quad
$d_{1} ^{(0)} \psi_{1 }^{(1)}(e_i^{(1)})=d_{1} ^{(0)}\left(x_{i-1}\varepsilon_0 + e_{i} \wedge (x_0^{\mu} e_{n-b} -x_n^a e_{n} ) \right)= $\\[2mm]
$ \left[ x_{i-1}g_0  +x_0^{\mu}(x_{i-1}x_{n-b}-x_i x_{n-b-1})-x_n^a(x_{i-1}x_{n}-x_i x_{n-1})\right]\varepsilon=x_ig_1 \varepsilon$.\\[2mm]
If  $s=2$, \ since \,  
 $d_2'(e_{i}^{(1)}\wedge e_{j}^{(1)})=\Delta_1(e_{i} \wedge e_{j} ) = x_{i}e_{j} ^{(1)}-x_{j}e_{i} ^{(1)} $\quad and (for $\,\, q=1,2)$ \\
$d_{2 }^{(0)}\left((e_i ^{(0)})\wedge e_j ^{(0)})\oplus (e_h\wedge e_k\wedge e_l \otimes  \lambda_q )\right)= 
(x_{i-1}x_{j}-x_{i}x_{j-1})\varepsilon_0 \oplus\left(- g_0 \ e_i\wedge e_j +d_2(e_h\wedge e_k\wedge e_l\otimes \lambda_q )\right),$  
\enu
\item[]$  \psi_1^{(1)}d_2'(e_{i}^{(1)}\wedge e_{j}^{(1)})=\psi_1^{(1)}(x_{i}e_{j}^{(1)}-x_{j}e_{i}^{(1)})=$\\
$
x_{i}\left[ x_{j-1}\varepsilon_0 +e_{j} \wedge (x_0^{\mu} e_{n-b} -x_n^a e_{n} ) \right] -x_{j}\left[ x_{i-1}\varepsilon_0 +e_{i} \wedge (x_0^{\mu} e_{n-b} -x_n^a e_{n} )\right]   =$\\
$ (x_{i}x_{j-1}-x_{i-1}x_{j} )\, \varepsilon_0  
 \oplus (x_{i} e_{j} -x_{j} e_{i}) \wedge (x_0^{\mu} e_{n-b} -x_n^a e_{n} );   $ 
\item[]$ d_{2 }^{(0)}\psi_{2}^{(1)}( e^{(1)}_{i}\wedge e^{(1)}_{j} )= $
$d_{2 }^{(0)}\left(-e_{i}^{(0)}\wedge e_{j}^{(0)}\oplus   e_i\wedge e_j\wedge  (x_0^{\mu}e_{n-b} 
-x_n^a  e_n )\otimes \lambda_1\right)= $\\
$-(x_{i-1}x_{j}-x_{i}x_{j-1})\, \varepsilon_0 \oplus g_0( e_i\wedge e_j )\, +\Delta_1\left(  e_i\wedge e_j\wedge  (x_0^{\mu}e_{n-b} 
-x_n^a  e_n ) \right).$\denu
 The conclusion and also the   commutativity for $s\geq 3$, 
   follow by the next Lemma \ref{delta}  (see the proof given in the following \lq\lq Step $h$'' which holds in the general case).\\[2mm]
Now, again by  the mapping cone construction, we get the exact complex:\\
\centerline
{$\begin{array}{lccccccccccccccccccccc}
	&&  d_{s} ^{(1)}&  &d_{1} ^{(1)} \\[-.5mm] 
	 {\mathbb F}_1:&    0\p   F^{(1)}_{n+1}\p...\p F^{(1)}_{s}  &\p    &  		F^{(1)}_{s-1}\p... \p   F^{(1)}_{1}& \p& F^{(1)}_0\p   P/  {\goth C} _1&\!\!	\longrightarrow\!\!& 0\\ 
\end{array} $}\\[2mm]
   where $\quad F^{(1)}_{s} =   K_{s-1}(-\delta_1)\oplus F^{(0)}_{s}= K_{s-1}(-\delta_1)\oplus E_{s-1}(-\delta_0)\oplus E_{s},\quad     {\goth C} _1=(  {\goth C} ,g_0,g_1),$  \\[2mm] 
   $  d_1 ^{(1)} :  K_{0}(-\delta_1)\oplus E_{0}(-\delta_0)\oplus E_{1}\p E_0(-\delta_0)$,\quad  
$ d_1 ^{(1)}( a_1\varepsilon_1\oplus a_0 \varepsilon_0 \oplus e_i\wedge e_j)=\left( a_1 g_1+a_0g_0+ x_{i-1}x_{j}-x_{i}x_{j-1}\right)\varepsilon_0 $\\[2mm]
$ d_s ^{(1)}= \left(
\begin{array}{ccc}
  d'_{s-1}&0       \\
  (-1)^{s-1}\psi^{(1)}_{s-1}& d^{(0)}_{s}  \\
     \end{array}
\right),\quad(s\geq 2)$. 
\begin{lemma}\label{delta}
Let $W:=e_{i_1} \wedge...\wedge e_{i_{s}},\quad  (i_1<i_2<...<i_s )$;  \, then for\,  $q\in\{0,1\}$, we have\quad 
 \enu
\item 
$\Delta_q(  W\wedge e_k)=\Delta_q(  W)\wedge e_k +(-1)^{s }x_{k+q-1}W$.
\item
 $d_s\circ\phi_s^{(h)}(W)=-\phi_{s-1}^{(h-1)}(\Delta_0 (W ))+\phi_{s-1}^{(h)}(\Delta_1 (W ))+   (-1)^{s+h}g_0 W \otimes \lambda_0^{h-1}\lambda_1^{s-h-1}.$
 \denu
 \end{lemma}
Proof.  
(1). \quad Let $p$ be the number of permutations to have $i_1<...i_{s-p}<i_k<i_{s-p+1}<...<i_s$. Then \\ $\Delta_q (  W\wedge e_k)=
(-1)^p[ \Delta_q (  e_{i_1} \wedge...e_{i_{s-p}}\wedge e_k\wedge  e_{i_{s-p+1}}\wedge...\wedge e_{i_{s}} )]=$\\[2mm] 
 $(-1)^p [\Delta_q (  e_{i_1} \wedge...  e_{i_{s-p}})\wedge e_k\wedge  e_{i_{s-p+1}}\wedge...\wedge e_{i_{s}} +(-1)^{s-p+2}x_{k+q-1}W+$\\$(-1)^{s-p+3} e_{i_1} \wedge...  e_{i_{s-p}}\wedge e_k\wedge\Delta_q (e_{i_{s-p+1}} \wedge... \wedge e_{i_{s}})]=$ \\[2mm]
 $ \Delta_q (  e_{i_1} \wedge...  e_{i_{s-p}})\wedge  e_{i_{s-p+1}}\wedge...\wedge e_{i_{s}}\wedge e_k +(-1)^{s+2}x_{k+q -1}W+(-1)^{s-p+2} e_{i_1} \wedge...  e_{i_{s-p}}\wedge\Delta_q (e_{i_{p+1}} \wedge...  e_{i_{s}})\wedge e_k=$\\[2mm] $\Delta_q (  W)\wedge e_k +(-1)^{s+2}x_{k+q -1}W.$\quad 
(Note that this result holds also for  $k\in\{i_1,...,i_s\}$).\\[2mm]
 (2).\quad Recall that 
  $\phi^h_s(W)=  \sum_{k=1}^h(-1)^k W\wedge  ( x_n^a e_{n+k-h} ^{}-x_0^{\mu} e_{n-b+k-h} ^{})\otimes \lambda_0^{k-1}\lambda_1^{s-k} $, hence\\[2mm]
$\begin{array}{lll} 
 d_s\circ\phi_s^{(h)}(W)=&\sum_{k=1}^h(-1)^k \Delta_0\left( W\wedge  ( x_n^a e_{n+k-h} ^{}-x_0^{\mu} e_{n-b+k-h} ^{})\right)\otimes \lambda_0^{k-2}\lambda_1^{s-k}+\\
 &\sum_{k=1}^h(-1)^k \Delta_1\left( W\wedge  ( x_n^a e_{n+k-h} ^{}-x_0^{\mu} e_{n-b+k-h} ^{})\right)\otimes \lambda_0^{k-1}\lambda_1^{s-k-1} 
 \end{array}$\\[2mm]
Let  $V_k:= x_n^a e_{n-h+k}  -x_0^{\mu} e_{n-h-b+k } $ and recall that
 $ g_h = x_n^ax_{n-h}-x_0^{\mu} x_{n-b-h}$;\ \ by  (\ref{delta}.1) we get\quad\\[2mm]     $\begin{array}{lll} d_s\circ\phi_s(W)=&
  \sum_{k=1}^h\left[(-1)^k \Delta_0( W)\wedge  V_k \otimes \lambda_0^{k-2}\lambda_1^{s-k}+(-1)^{s }g_{h-k+1} W\otimes \lambda_0^{k-2}\lambda_1^{s-k}\right]+\\
  &\sum_{k=1}^h\left[(-1)^k \Delta_1 ( W)\wedge  V_k \otimes \lambda_0^{k-1}\lambda_1^{s-k-1}+(-1)^{s } g_{h-k}    W\otimes \lambda_0^{k-1}\lambda_1^{s-k-1}\right]\ = \end{array}$\\[2mm]
  $\begin{array}{lll} 
\sum_{k=2}^h(-1)^k \Delta_0( W)\wedge  V_k \otimes \lambda_0^{k-2}\lambda_1^{s-k}+\phi_{s-1}^{h}(\Delta_1 (W))+\\
  
  (-1)^{s} \sum_{k=1}^h\ (-1)^{k}\left[(g_{h-k+1} )W \otimes \lambda_0^{k-2}\lambda_1^{s-k} +  (g_{h-k}  )W \otimes \lambda_0^{k-1}\lambda_1^{s-k-1} \right] = \end{array}$\\[2mm]
  $\begin{array}{lll} 
 \sum_{k'=1}^{h-1}(-1)^{k'+1} \Delta_0( W)\wedge  ( x_n^a e_{n+k'+1-h}  -x_0^{\mu} e_{n-b+k'+1-h} ) \otimes \lambda_0^{k'-1}\lambda_1^{s-k'-1}+\phi_{s-1}^{h}(\Delta_1 (W))+\\
  
  (-1)^{s}[ -   (x_n^a x_{n+1-h} -x_0^{\mu} x_{n-b+1-h} )W \otimes  \lambda_1^{s-2}+\\
  (x_n^a x_{n+1-h } -x_0^{\mu} x_{n-b+1-h } ) W\otimes  \lambda_1^{s-2} +  (x_n^a x_{n+2-h} -x_0^{\mu} x_{n-b+2-h} ) W\otimes \lambda_0 \lambda_1^{s-3} +\\
 -(x_n^a x_{n+2-h } -x_0^{\mu} x_{n-b+2-h} ) W\otimes \lambda_0\lambda_1^{s-3} +  (x_n^a x_{n+3-h} -x_0^{\mu} x_{n-b+3-h} ) W\otimes \lambda_0^{2}\lambda_1^{s-4}  +...\\
 (-1)^{h}(x_n^a x_{n-1} -x_0^{\mu} x_{n-b-1} )W \otimes \lambda_0^{h-2}\lambda_1^{s-h} +(-1)^{h}  (x_n^{a+1} -x_0^{\mu} x_{n-b} )W \otimes \lambda_0^{h-1}\lambda_1^{s-h-1}] = \end{array}$\\[2mm]
 
 $\begin{array}{lll} 
=-\phi_{s-1}^{(h-1)}(  \Delta_0(W))+\phi_{s-1}^{h}(\Delta_1 (W))+   (-1)^{s+h}g_0 W \otimes \lambda_0^{h-1}\lambda_1^{s-h-1}\ \qquad ({\rm where \,\, }k':=k-1). \quad\diamond\end{array}$ \\[2mm]

{\bf   STEP \, h  \, ($2\leq h\leq n-b$).}
 \, By iterating this method, for all $h=1,...,n-b $, let   $  {\goth C} _{h}=  {\goth C} +(g_0,...,g_{h})$. \\
 By means of the following   commutative diagram\\
  \centerline{$\begin{array}{cccccccccccccccccccccc}
 &&& &&&&& 0\\
& &&  &&&& &\downarrow \\
& && d'_{s}&& &&  d'_0\\  
 {\mathbb K}(-\delta_{h}) :&...&      K_{s}(-\delta_{h})&\!\!\longrightarrow\!\!& K_{s-1}(-\delta_{h})&...    &K_0(-\delta_{h})&\!\!\longrightarrow\!\!& P(-\delta_{h})/(  {\goth C} _{h-1}:g_h)&\!\!\longrightarrow\!\!& 0\\
\\
 & &\downarrow\psi_{s }^{(h)}&&\downarrow\psi_{s-1}^{(h)}&  &\downarrow\psi_{0}^{(h)}&&\ \ \ \ \downarrow \cdot g_h \\
 & &&  d_{s} ^{(h-1)}& & &&d_0 ^{(h-1)}  \\ 
 {\mathbb F}_{h-1}:&  ...&   F^{(h-1)}_{s}  \ \ \    &\!\!\longrightarrow\!\!& F^{(h-1)}_{s-1}\ \ \ \ \   &...  & F^{(h-1)}_0  \ \ \ \ \ &\!\!\longrightarrow\!\!& P/  {\goth C} _{h-1}&\!\!\longrightarrow\!\!& 0\\
\\ &&& & &&& \searrow&\downarrow \\ \\
  &&&&&&& &P/  {\goth C} _h \\ \\
  & &&&&& &&\downarrow \\
&&&&&& &&0
  \end{array} $}\\[2mm]
   we construct the free complex ${\mathbb F}_h $:\\ 
  $\begin{array}{cccccccccccccccccccccc}

 & && & d_{s} ^{(h )}& & &&&&d_0 ^{(h )}  \\ 

 {\mathbb F}_{h }:\quad 0\p F^{(h)}_{n}    &\p&  ...&   F^{(h)}_{s}    &\p& F^{(h)}_{s-1}&\p &... &\p & F^{(h)}_0  &\p& P/  {\goth C} _{h}&\!\!\longrightarrow\!\!& 0\\
  \end{array} $\\[2mm]
  with   \centerline{$\begin{array}{llllllllll}
    F_0^{(h)}=E_0(-\delta_h),\ &&
    F_s^{(h )}=K_{s-1}(-\delta_{h })\oplus...\oplus K_{s-1}(-\delta_{1})\oplus E_{s-1}(-\delta_0)\oplus E_s \quad (s,h\geq 1)\\
  \end{array}$.}\\[2mm]
 The commutativity for $s\geq 2$ of the above diagram follows by   Lemma \ref{delta},  since, with notation (\ref{nota1}.6)
 \enu
 
\item[A.] 

$\psi_{s-1}^{(h)}d'_s(W)=
 \psi_{s-1}^{(h)}\Delta_1(W)=
\left[\begin{array}{lllllcccccccc}

&(-1)^{s}\Delta_0(\Delta_1( W))&  &[\in K_{s-2}(-\delta_{h-1}) ]\\

 \oplus&(-1)^{h}  \Delta_1( W)\otimes\lambda_0^{h-1}\lambda_1^{s-h-2}& &[\in E_{s-2}(-\delta_0)] \\
 \oplus&  \phi^{(h)}_{s-1}(\Delta_1(W))&&[\in   E_{s-1}]
 
  \end{array}\right. $

\item[B.] $d_s^{(h-1)}\psi^{(h)}_{s}(W)=d_s^{(h-1)} \left[  (-1)^{s-1}\Delta_0(W)\oplus (-1)^hW\otimes\lambda_0^{h-1}\lambda_1^{s-h-1} \oplus \phi^{(h)}_s (W ) \right] =$

\item[] $\Delta_1(  (-1)^{s-1}\Delta_0(W))\oplus \Delta_0( (-1)^{s }\Delta_0(W)) $\\[2mm]
$
\oplus (-1)^{s+h}    (-1)^{s-1}\Delta_0(W) \otimes\lambda_0^{h-2}\lambda_1^{s-h-1 } + d_{s-1}((-1)^hW\otimes\lambda_0^{h-1}\lambda_1^{s-h-1})  $\\[2mm]\
  $
 \oplus (-1)^{s-1}\phi_{s-1}^{(h-1)}( (-1)^{s-1}\Delta_0(W))+ (-1)^{s-1}(-1)^hW\otimes\lambda_0^{h-1}\lambda_1^{s-h-1}\, g_0 + d_s(\phi^{(h)}_s (W )) .$

\item[]$d_{s}^{(h-1)}\psi_{s}^{(h)}(W)=
\left[\begin{array}{lllllcccccccc}
&\Delta_1(  (-1)^{s-1}\Delta_0(W))\\

 \oplus &    (-1)^{s }\Delta^2_0(W))  \\
 
  \oplus &   (-1)^{h-1}     \Delta_0(W) \otimes\lambda_0^{h-2}\lambda_1^{s-h-1 } + d_{s-1}((-1)^hW\otimes\lambda_0^{h-1}\lambda_1^{s-h-1})  \\
  
  \oplus&  \phi_{s-1}^{(h-1)}(  \Delta_0(W))+ (-1)^{s+h-1} W\otimes\lambda_0^{h-1}\lambda_1^{s-h-1}\, g_0 + d_s(\phi^{(h)}_s (W )).
 \quad\diamond\end{array}\right.$

\denu
\begin{rem}\label{redund}{\rm
The  complex ${\mathbb F}_0$ itself is  a minimal free resolution for the ideal $I$ if $\, b=n $, but not in the other cases.\quad
In fact, if $\,b<n$ \, and $s\geq 2$  the maps $  \psi_{s}^{(h)} $  give invertible entries in the matrix $M_{s+1}$ associated to the map $d_{s+1}^{(n-b)}$ (note that in $M_2$ all the entries are non-invertible according to the definition of $\psi_{1}^{(1)}$).\\ 
To obtain a minimal free resolution we have to delete suitable subpaces of the modules $C_s$   (see \ref{cs}):
the invertible  entries in $\, M_{s+1}$ arise as $\psi_{s}^{(h)}(W^{(h)})\cap E_{s-1}(-\delta_0)$ \quad i.e.\, 
  $(-1)^h W^{(0)}\otimes\lambda_0^{h-1}\lambda_1^{s-h-1}, h=1,...,n-b$. Since such elements are considered only with non-negative   powers of $\lambda_0,\lambda_1$, to compute the dimension of the redundant subspace $D_s$ define\\ 
  \centerline{$\nu_s:= min \{s-1, n-b\}$. } Then\ \ $ \nu_s\binom{n}{s}=dim\left(Im\, \psi_{s}^{(h)} \cap E_{s-1}(-\delta_0)\right) $.\ \
  Further $\psi_{s-1}^{(h)} $ gives invertible entries in $M_s$.\\
  Therefore for each $s\geq 2$, the redundant subspace $D_s\subseteq C_s$ has dimension $  \nu_s\binom{n}{s}+ \nu_{s-1}\binom{n}{s-1},\quad (s\geq 2)$:}
  \end{rem}
    \begin{thm}
   A minimal free resolution of the ideal I defining an $AS$ curve with $b>1$  is given by the complex
   \\
   \centerline{$\begin{array}{cccccccccccccccccccccc}
 &&   &&&  &&  &&  \\ 
 { \cal R} :&&  0 \p    R_{n}  \p ...\p R_{s} &\p... \p & R_2  &\!\!\longrightarrow\!\!& R_1  &\!\!\longrightarrow\!\!& P\,    \end{array} $}\\[2mm]
where \,\,$R_{s}=C_s/D_s$ \,\ according to Remark \ref{redund}; \quad in particular\\[2mm]
$\left[\begin{array}{llll}
R_n &=&\bigoplus_{0\leq k\leq b-2 }\left( e_1\wedge...\wedge e_n\, \lambda_0^{n-b+k}\lambda_1^{b-k-2}\right)P,  \quad {\rm if}\,\, b\geq 2\\
&=& K_{n-1}(-\delta_1),\quad {\rm if}\,\, b=1\\[2mm]
R_2&=& K_{1}(-\delta_{n-b })\oplus...\oplus K_{1}(-\delta_{1})\oplus E_2, \quad {\rm with}\\
  &&K_{1}(-\delta_{j })=\bigoplus_{1\leq i\leq n } P\left(-m_i-\delta_j\right),\\
  &&E_2=\bigoplus_{1\leq i_1<i_2<i_3\leq n, q=0,1} P\big(-m_{i_1}-m_{i_2}-m_{i_3}+(2-q)d\, \big)\\[2mm] 
R_1&=& K_{0}(-\delta_{n-b })\oplus...\oplus K_{0}(-\delta_{1})\oplus E_{0}(-\delta_0)\oplus E_1 \\
&\simeq& \bigoplus_{0\leq i\leq n-b } P\left(-\delta_i\right ) \bigoplus_{1\leq i_1<i_2 \leq n } P\left(-m_{i_1}-m_{i_2}+d \right )\,.
\\

\end{array} \right.$

\end{thm}
\begin{coro}\label{syz} \enu \item The Betti numbers of the ideal $I$ are \\[2mm]
\centerline{$\beta_s=dim(R_s)= \left[\begin{array}{llll}
(n-b+2-s)\binom{n}{s-1}+ s\binom{n}{s+1} , & if  & 2\leq s< n-b+ 2   \\
\\
 (s-1-n+b)\binom{n}{s }+ s\binom{\,\, n\,\,}{s+1} ,& if  &  n-b+2\leq s\leq n\,.
 \end{array}\right.$}
 \item The local ring $R=k[[x_0,...,x_n]]/I$ where $I$ is the defining ideal of an $AS$ monomial curve is of homogeneous type, i.e., \, $\beta_i(R)=\beta_i\big({\cal G}\big)$, where   ${\cal G}$ is the associated graded ring of $R$. 
 \denu
 \end{coro}
Proof. (1). It suffices to recall that\\
 $dim(R_s)= dim(C_s)-dim (D_s)=(n-b)\binom{n}{s-1}+(s-1)\binom{n}{s}+ s\binom{n}{s+1}- \nu_s\binom{n}{s}- \nu_{s-1}\binom{n}{s-1},\quad (s\geq 2)$.\\
 Since $\nu_s:= min \{s-1, n-b\}$,  we get\\ \centerline{$\beta_s=dim(R_s)= (n-b- \nu_{s-1})\binom{n}{s-1}+(s-1-\nu_s)\binom{n}{s}+ s\binom{n}{s+1}.$} \\[2mm]
 (2)  follows immediately by \cite[Theorem 4.1]{lr}
  \quad$\diamond$

\renewcommand{\refname}{{\large\bf References}}

\end{document}